\documentclass{amsproc}
\usepackage{amssymb}
\usepackage{amsmath}
\usepackage{eurosym}
\usepackage{amsfonts}

            \advance \topmargin by -\headheight \advance
            \topmargin by -\headsep
            \evensidemargin
            \oddsidemargin
\newtheorem{theorem}{Theorem}[section]
\newtheorem{lemma}[theorem]{Lemma}
\newtheorem{definition}[theorem]{Definition}
\newtheorem{corollary}[theorem]{Corollary}

\numberwithin{equation}{section}
\input{tcilatex}

\begin{document}
\title[Matrix Sequences of Tribonacci and Tribonacci-Lucas Numbers]{Matrix
Sequences of Tribonacci and Tribonacci-Lucas Numbers}
\thanks{}
\author[Y\"{u}ksel ~Soykan]{Y\"{U}KSEL SOYKAN}
\maketitle

\begin{center}
\textsl{Zonguldak B\"{u}lent Ecevit University, Department of Mathematics, }

\textsl{Art and Science Faculty, 67100, Zonguldak, Turkey }

\textsl{e-mail: \ yuksel\_soykan@hotmail.com}
\end{center}

\textbf{Abstract.} In this paper, we define Tribonacci and Tribonacci-Lucas
matrix sequences and investigate their properties.



\section{Introduction and Preliminaries}

Recently, there have been so many studies of the sequences of numbers in the
literature that concern about subsequences of the Horadam numbers and
generalized Tribonacci numbers such as Fibonacci, Lucas, Pell and Jacobsthal
numbers; Tribonacci, Tribonacci-Lucas, Padovan, Perrin, Padovan-Perrin,
Narayana, third order Jacobsthal and third order Jacobsthal-Lucas numbers.

The sequences of numbers were widely used in many research areas, such as
physics, engineering, architecture, nature and art. For example, the ratio
of two consecutive Fibonacci numbers converges to the Golden section
(ratio), $\alpha _{F}=\frac{1+\sqrt{5}}{2}$; which appears in modern
research, particularly physics of the high energy particles or theoretical
physics. Another example, the ratio of two consecutive Padovan numbers
converges to the Plastic ratio, $\alpha _{P}=\sqrt[3]{\frac{1}{2}+\frac{1}{6}%
\sqrt{\frac{23}{3}}}+\sqrt[3]{\frac{1}{2}-\frac{1}{6}\sqrt{\frac{23}{3}}},$
which have many applications to such as architecture, see [\ref%
{bib:marohnic2012}]. One last example, the ratio of two consecutive
Tribonacci numbers converges to the Tribonacci ratio, $\alpha _{T}=\frac{1+%
\sqrt[3]{19+3\sqrt{33}}+\sqrt[3]{19-3\sqrt{33}}}{3}.$ For a short
introduction to these three constants, see [\ref{piezas}].

On the other hand, the matrix sequences have taken so much interest for
different type of numbers. For matrix sequences of generalized Horadam type
numbers, see for example [\ref{civciv2008}], [\ref{civciv2008b}], [\ref%
{gulec2012}], [\ref{uslu2013}], [\ref{uygun2015}], [\ref{uygun2016}], [\ref%
{yazlik2013a}], [\ref{wani2018}], and for matrix sequences of generalized
Tribonacci type numbers, see for instance [\ref{cerdamoralez2018aa}], [\ref%
{yilmaz2013}], [\ref{yilmaz2014a}].

In this paper, the matrix sequences of Tribonacci and Tribonacci-Lucas
numbers will be defined for the first time in the literature. Then, by
giving the generating functions, the Binet formulas, and summation formulas
over these new matrix sequences, we will obtain some fundamental properties
on Tribonacci and Tribonacci-Lucas numbers. Also, we will present the
relationship between these matrix sequences.

First, we give some background about Tribonacci and Tribonacci-Lucas
numbers. Tribonacci sequence $\{T_{n}\}_{n\geq 0}$ (sequence A000073 in [\ref%
{bib:sloane}]) and Tribonacci-Lucas sequence $\{K_{n}\}_{n\geq 0}$ (sequence
A001644 in [\ref{bib:sloane}]) are defined by the third-order recurrence
relations 
\begin{equation}
T_{n}=T_{n-1}+T_{n-2}+T_{n-3},\text{ \ \ \ \ }T_{0}=0,T_{1}=1,T_{2}=1,
\label{equati:fvcvxghsbnz}
\end{equation}%
and 
\begin{equation}
K_{n}=K_{n-1}+K_{n-2}+K_{n-3},\text{ \ \ \ \ }K_{0}=3,K_{1}=1,K_{2}=3
\label{equati:pazertvbcunsmn}
\end{equation}%
respectively. Tribonacci concept was introduced by M. Feinberg\ [\ref%
{bib:feinberg1963}] in 1963. Basic properties of it is given in [\ref%
{bib:bruce1984}], [\ref{bib:scott1977}], [\ref{bib:shannon1977}], [\ref%
{bib:yalavigi1972}] and Binet formula for the $n$th number is given in [\ref%
{bib:spickerman1981}].

The sequences $\{T_{n}\}_{n\geq 0}$ and $\{K_{n}\}_{n\geq 0}$ can be
extended to negative subscripts by defining%
\begin{equation*}
T_{-n}=-T_{-(n-1)}-T_{-(n-2)}+T_{-(n-3)}
\end{equation*}%
and%
\begin{equation*}
K_{-n}=-K_{-(n-1)}-K_{-(n-2)}+K_{-(n-3)}
\end{equation*}%
for $n=1,2,3,...$ respectively. Therefore, recurrences (\ref%
{equati:fvcvxghsbnz}) and (\ref{equati:pazertvbcunsmn}) hold for all integer 
$n.$

By writting $T_{n-1}=T_{n-2}+T_{n-3}+T_{n-4}$ and eliminating $T_{n-2}$ and $%
T_{n-3}$ between this recurrence relation and recurrence relation (\ref%
{equati:fvcvxghsbnz}), a useful alternative recurrence relation is obtained
for $n\geq 4:$%
\begin{equation}
T_{n}=2T_{n-1}-T_{n-4},\ \ \ \ T_{0}=0,T_{1}=T_{2}=1,T_{3}=2.
\label{equati:dostyarewqdcxz}
\end{equation}%
Extension of the definition of $T_{n}$ to negative subscripts can be proved
by writing the recurrence relation (\ref{equati:dostyarewqdcxz}) as%
\begin{equation*}
T_{-n}=2T_{-n+3}-T_{-n+4}.
\end{equation*}%
Note that $T_{-n}=T_{n-1}^{2}-T_{n-2}T_{n},$ (see [\ref{bib:choi2013}]).

We can give some relations between $\{T_{n}\}$ and $\{K_{n}\}$\ as 
\begin{equation}
K_{n}=3T_{n+1}-2T_{n}-T_{n-1}  \label{equation:gfdvbxczvsadou}
\end{equation}%
and%
\begin{equation}
K_{n}=T_{n}+2T_{n-1}+3T_{n-2}  \label{equation:fvcdxsartewqa}
\end{equation}%
and also%
\begin{equation}
K_{n}=4T_{n+1}-T_{n}-T_{n+2}.  \label{equation:yuhtgsdafcscvb}
\end{equation}%
Note that the last three identities hold for all integers $n.$

The first few Tribonacci numbers and Tribonacci Lucas numbers with positive
subscript are given in the following table:

$%
\begin{array}{ccccccccccccccc}
n & 0 & 1 & 2 & 3 & 4 & 5 & 6 & 7 & 8 & 9 & 10 & 11 & 12 & ... \\ 
T_{n} & 0 & 1 & 1 & 2 & 4 & 7 & 13 & 24 & 44 & 81 & 149 & 274 & 504 & ... \\ 
T_{-n} & 0 & 0 & 1 & -1 & 0 & 2 & -3 & 1 & 4 & -8 & 5 & 7 & -20 & ...%
\end{array}%
$

The first few Tribonacci numbers and Tribonacci Lucas numbers with negative
subscript are given in the following table:

$%
\begin{array}{ccccccccccccccc}
n & 0 & 1 & 2 & 3 & 4 & 5 & 6 & 7 & 8 & 9 & 10 & 11 & 12 & ... \\ 
K_{n} & 3 & 1 & 3 & 7 & 11 & 21 & 39 & 71 & 131 & 241 & 443 & 815 & 1499 & 
... \\ 
K_{-n} & 3 & -1 & -1 & 5 & -5 & -1 & 11 & -15 & 3 & 23 & -41 & 21 & 43 & ...%
\end{array}%
$

It is well known that for all integers $n,$ usual Tribonacci and
Tribonacci-Lucas numbers can be expressed using Binet's formulas%
\begin{equation}
T_{n}=\frac{\alpha ^{n+1}}{(\alpha -\beta )(\alpha -\gamma )}+\frac{\beta
^{n+1}}{(\beta -\alpha )(\beta -\gamma )}+\frac{\gamma ^{n+1}}{(\gamma
-\alpha )(\gamma -\beta )}  \label{equat:mnopcvbedcxzsa}
\end{equation}%
and%
\begin{equation}
K_{n}=\alpha ^{n}+\beta ^{n}+\gamma ^{n}  \label{equation:cfrdcsxszouea}
\end{equation}%
respectively, where $\alpha ,\beta $ and $\gamma $ are the roots of the
cubic equation $x^{3}-x^{2}-x-1=0.$ Moreover, 
\begin{eqnarray*}
\alpha &=&\frac{1+\sqrt[3]{19+3\sqrt{33}}+\sqrt[3]{19-3\sqrt{33}}}{3}, \\
\beta &=&\frac{1+\omega \sqrt[3]{19+3\sqrt{33}}+\omega ^{2}\sqrt[3]{19-3%
\sqrt{33}}}{3}, \\
\gamma &=&\frac{1+\omega ^{2}\sqrt[3]{19+3\sqrt{33}}+\omega \sqrt[3]{19-3%
\sqrt{33}}}{3}
\end{eqnarray*}%
where%
\begin{equation*}
\omega =\frac{-1+i\sqrt{3}}{2}=\exp (2\pi i/3),
\end{equation*}%
is a primitive cube root of unity. Note that we have the following identities%
\begin{eqnarray*}
\alpha +\beta +\gamma &=&1, \\
\alpha \beta +\alpha \gamma +\beta \gamma &=&-1, \\
\alpha \beta \gamma &=&1.
\end{eqnarray*}

The generating functions for the Tribonacci sequence $\{T_{n}\}_{n\geq 0}$
and Tribonacci-Lucas sequence $\{K_{n}\}_{n\geq 0}$ are%
\begin{equation}
\sum_{n=0}^{\infty }T_{n}x^{n}=\frac{x}{1-x-x^{2}-x^{3}}\text{ \ and \ }%
\sum_{n=0}^{\infty }K_{n}x^{n}=\frac{3-2x-x^{2}}{1-x-x^{2}-x^{3}}.
\label{equation:yugdfvxbgsopqac}
\end{equation}

Note that the Binet form of a sequence satisfying (\ref{equati:fvcvxghsbnz}%
)\ and (\ref{equati:pazertvbcunsmn}) for non-negative integers is valid for
all integers $n.$ This result of Howard and Saidak [\ref{bib:howard2010}] is
even true in the case of higher-order recurrence relations as the following
theorem shows.

\begin{theorem}
\label{theorem:fvgxdfsxczsaer}[\ref{bib:howard2010}]Let $\{w_{n}\}$ be a
sequence such that 
\begin{equation*}
\{w_{n}\}=a_{1}w_{n-1}+a_{2}w_{n-2}+...+a_{k}w_{n-k}
\end{equation*}%
for all integers $n,$ with arbitrary initial conditions $%
w_{0},w_{1},...,w_{k-1}.$ Assume that each $a_{i}$ and the initial
conditions are complex numbers. Write%
\begin{eqnarray}
f(x) &=&x^{k}-a_{1}x^{k-1}-a_{2}x^{k-2}-...-a_{k-1}x-a_{k}
\label{equation:mnbvyuhgoewapvbc} \\
&=&(x-\alpha _{1})^{d_{1}}(x-\alpha _{2})^{d_{2}}...(x-\alpha _{h})^{d_{h}} 
\notag
\end{eqnarray}%
with $d_{1}+d_{2}+...+d_{h}=k,$ and $\alpha _{1},\alpha _{2},...,\alpha _{k}$
distinct. Then

\begin{description}
\item[(a)] For all $n,$%
\begin{equation}
w_{n}=\sum_{m=1}^{k}N(n,m)(\alpha _{m})^{n}
\label{equation:yuhnbvdfscxzratqw}
\end{equation}%
where%
\begin{equation*}
N(n,m)=A_{1}^{(m)}+A_{2}^{(m)}n+...+A_{r_{m}}^{(m)}n^{r_{m}-1}=%
\sum_{u=0}^{r_{m}-1}A_{u+1}^{(m)}n^{u}
\end{equation*}%
with each $A_{i}^{(m)}$ a constant determined by the initial conditions for $%
\{w_{n}\}$. Here, equation (\ref{equation:yuhnbvdfscxzratqw}) is called the
Binet form (or Binet formula) for $\{w_{n}\}.$ We assume that $f(0)\neq 0$
so that $\{w_{n}\}$ can be extended to negative integers $n.$

If the zeros of (\ref{equation:mnbvyuhgoewapvbc}) are distinct, as they are
in our examples, then%
\begin{equation*}
w_{n}=A_{1}(\alpha _{1})^{n}+A_{2}(\alpha _{2})^{n}+...+A_{k}(\alpha
_{k})^{n}.
\end{equation*}

\item[(b)] The Binet form for $\{w_{n}\}$ is valid for all integers $n.$
\end{description}
\end{theorem}

\section{The Matrix Sequences of Tribonacci and Tribonacci-Lucas Numbers}

In this section we define Tribonacci and Tribonacci-Lucas matrix sequences
and investgate their properties.

\begin{definition}
\label{definition:tyfsdczxawuqp}For any integer $n\geq 0,$ the Tribonacci
matrix $(\mathcal{T}_{n})$ and Tribonacci-Lucas matrix $(\mathcal{K}_{n})$
are defined by%
\begin{eqnarray}
\mathcal{T}_{n} &=&\mathcal{T}_{n-1}+\mathcal{T}_{n-2}+\mathcal{T}_{n-3},
\label{equation:yusoabpsmnbscv} \\
\mathcal{K}_{n} &=&\mathcal{K}_{n-1}+\mathcal{K}_{n-2}+\mathcal{K}_{n-3,}
\label{equation:dfscxzvayuewrtsfg}
\end{eqnarray}%
respectively, with initial conditions%
\begin{equation*}
\mathcal{T}_{0}=\left( 
\begin{array}{ccc}
1 & 0 & 0 \\ 
0 & 1 & 0 \\ 
0 & 0 & 1%
\end{array}%
\right) ,\mathcal{T}_{1}=\left( 
\begin{array}{ccc}
1 & 1 & 1 \\ 
1 & 0 & 0 \\ 
0 & 1 & 0%
\end{array}%
\right) ,\mathcal{T}_{2}=\left( 
\begin{array}{ccc}
2 & 2 & 1 \\ 
1 & 1 & 1 \\ 
1 & 0 & 0%
\end{array}%
\right)
\end{equation*}%
and%
\begin{equation*}
\mathcal{K}_{0}=\left( 
\begin{array}{ccc}
1 & 2 & 3 \\ 
3 & -2 & -1 \\ 
-1 & 4 & -1%
\end{array}%
\right) ,\mathcal{K}_{1}=\left( 
\begin{array}{ccc}
3 & 4 & 1 \\ 
1 & 2 & 3 \\ 
3 & -2 & -1%
\end{array}%
\right) ,\mathcal{K}_{2}=\left( 
\begin{array}{ccc}
7 & 4 & 3 \\ 
3 & 4 & 1 \\ 
1 & 2 & 3%
\end{array}%
\right) .
\end{equation*}
\end{definition}

The sequences $\{\mathcal{T}_{n}\}_{n\geq 0}$ and $\{\mathcal{K}%
_{n}\}_{n\geq 0}$ can be extended to negative subscripts by defining 
\begin{equation*}
\mathcal{T}_{-n}=-\mathcal{T}_{-(n-1)}-\mathcal{T}_{-(n-2)}+\mathcal{T}%
_{-(n-3)}
\end{equation*}%
and%
\begin{equation*}
\mathcal{K}_{-n}=-\mathcal{K}_{-(n-1)}-\mathcal{K}_{-(n-2)}+\mathcal{K}%
_{-(n-3)}
\end{equation*}%
for $n=1,2,3,...$ respectively. Therefore, recurrences (\ref%
{equation:yusoabpsmnbscv}) and (\ref{equation:dfscxzvayuewrtsfg}) hold for
all integers $n.$

The following theorem gives the $n$th general terms of the Tribonacci and
Tribonacci-Lucas matrix sequences.

\begin{theorem}
\label{theorem:hbnvxgfdtreopas}For any integer $n\geq 0,$ we have the
following formulas of the matrix sequences:%
\begin{eqnarray}
\mathcal{T}_{n} &=&\left( 
\begin{array}{ccc}
T_{n+1} & T_{n}+T_{n-1} & T_{n} \\ 
T_{n} & T_{n-1}+T_{n-2} & T_{n-1} \\ 
T_{n-1} & T_{n-2}+T_{n-3} & T_{n-2}%
\end{array}%
\right)  \label{equat:yaczxdsweauyhnbv} \\
\mathcal{K}_{n} &=&\left( 
\begin{array}{ccc}
K_{n+1} & K_{n}+K_{n-1} & K_{n} \\ 
K_{n} & K_{n-1}+K_{n-2} & K_{n-1} \\ 
K_{n-1} & K_{n-2}+K_{n-3} & K_{n-2}%
\end{array}%
\right) .  \label{equati:bndexvsuzosavrw}
\end{eqnarray}
\end{theorem}

Proof. We prove (\ref{equat:yaczxdsweauyhnbv}) by strong mathematical
induction on $n$. (\ref{equati:bndexvsuzosavrw}) can be proved similarly.

If $n=0$ then, since $T_{1}=1,T_{2}=1,T_{0}=T_{-1}=0$,\ $T_{-2}=1,$ $%
T_{-3}=-1,$ we have 
\begin{equation*}
\mathcal{T}_{0}=\left( 
\begin{array}{ccc}
T_{1} & T_{0}+T_{-1} & T_{0} \\ 
T_{0} & T_{-1}+T_{-2} & T_{-1} \\ 
T_{-1} & T_{-2}+T_{-3} & T_{-2}%
\end{array}%
\right) =\left( 
\begin{array}{ccc}
1 & 0 & 0 \\ 
0 & 1 & 0 \\ 
0 & 0 & 1%
\end{array}%
\right)
\end{equation*}%
which is true and 
\begin{equation*}
\mathcal{T}_{1}=\left( 
\begin{array}{ccc}
T_{2} & T_{1}+T_{0} & T_{1} \\ 
T_{1} & T_{0}+T_{-1} & T_{0} \\ 
T_{0} & T_{-1}+T_{-2} & T_{-1}%
\end{array}%
\right) =\left( 
\begin{array}{ccc}
1 & 1 & 1 \\ 
1 & 0 & 0 \\ 
0 & 1 & 0%
\end{array}%
\right)
\end{equation*}%
which is true. Assume that the equality holds for $n\leq k.$ For $n=k+1,$ we
have%
\begin{eqnarray*}
\mathcal{T}_{k+1} &=&\mathcal{T}_{k}+\mathcal{T}_{k-1}+\mathcal{T}_{k-2} \\
&=&\left( 
\begin{array}{ccc}
T_{k+1} & T_{k}+T_{k-1} & T_{k} \\ 
T_{k} & T_{k-1}+T_{k-2} & T_{k-1} \\ 
T_{k-1} & T_{k-2}+T_{k-3} & T_{k-2}%
\end{array}%
\right) +\left( 
\begin{array}{ccc}
T_{k} & T_{k-1}+T_{k-2} & T_{k-1} \\ 
T_{k-1} & T_{k-2}+T_{k-3} & T_{k-2} \\ 
T_{k-2} & T_{k-3}+T_{k-4} & T_{k-3}%
\end{array}%
\right) \\
&&+\left( 
\begin{array}{ccc}
T_{k-1} & T_{k-2}+T_{k-3} & T_{k-2} \\ 
T_{k-2} & T_{k-3}+T_{k-4} & T_{k-3} \\ 
T_{k-3} & T_{k-4}+T_{k-5} & T_{k-4}%
\end{array}%
\right) \\
&=&\left( 
\begin{array}{ccc}
T_{k}+T_{k-1}+T_{k+1} & T_{k}+T_{k-1}+T_{k-1}+T_{k-2}+T_{k-2}+T_{k-3} & 
T_{k}+T_{k-1}+T_{k-2} \\ 
T_{k}+T_{k-1}+T_{k-2} & T_{k-1}+T_{k-2}+T_{k-2}+T_{k-3}+T_{k-3}+T_{k-4} & 
T_{k-1}+T_{k-2}+T_{k-3} \\ 
T_{k-1}+T_{k-2}+T_{k-3} & T_{k-2}+T_{k-3}+T_{k-3}+T_{k-4}+T_{k-4}+T_{k-5} & 
T_{k-2}+T_{k-3}+T_{k-4}%
\end{array}%
\right) \\
&=&\left( 
\begin{array}{ccc}
T_{k+2} & T_{k}+T_{k+1} & T_{k+1} \\ 
T_{k+1} & T_{k}+T_{k-1} & T_{k} \\ 
T_{k} & T_{k-1}+T_{k-2} & T_{k-1}%
\end{array}%
\right) .
\end{eqnarray*}%
Thus, by strong induction on $n,$ this proves (\ref{equat:yaczxdsweauyhnbv}%
). 
\endproof%

We now give the Binet formulas for the Tribonacci and Tribonacci-Lucas
matrix sequences.

\begin{theorem}
\label{theorem:sacasdmnbhgvc}For every integer $n,$ the Binet formulas of
the Tribonacci and Tribonacci-Lucas matrix sequences are given by%
\begin{eqnarray}
\mathcal{T}_{n} &=&A_{1}\alpha ^{n}+B_{1}\beta ^{n}+C_{1}\gamma ^{n},
\label{equation:mousxzadfgytc} \\
\mathcal{K}_{n} &=&A_{2}\alpha ^{n}+B_{2}\beta ^{n}+C_{2}\gamma ^{n}.
\label{equat:fvbxdszuqwsazx}
\end{eqnarray}%
where%
\begin{eqnarray*}
A_{1} &=&\frac{\alpha \mathcal{T}_{2}+\alpha (\alpha -1)\mathcal{T}_{1}+%
\mathcal{T}_{0}}{\alpha \left( \alpha -\gamma \right) \left( \alpha -\beta
\right) },B_{1}=\frac{\beta \mathcal{T}_{2}+\beta (\beta -1)\mathcal{T}_{1}+%
\mathcal{T}_{0}}{\beta \left( \beta -\gamma \right) \left( \beta -\alpha
\right) },C_{1}=\frac{\gamma \mathcal{T}_{2}+\gamma (\gamma -1)\mathcal{T}%
_{1}+\mathcal{T}_{0}}{\gamma \left( \gamma -\beta \right) \left( \gamma
-\alpha \right) } \\
A_{2} &=&\frac{\alpha \mathcal{K}_{2}+\alpha (\alpha -1)\mathcal{K}_{1}+%
\mathcal{K}_{0}}{\alpha \left( \alpha -\gamma \right) \left( \alpha -\beta
\right) },B_{2}=\frac{\beta \mathcal{K}_{2}+\beta (\beta -1)\mathcal{K}_{1}+%
\mathcal{K}_{0}}{\beta \left( \beta -\gamma \right) \left( \beta -\alpha
\right) },C_{2}=\frac{\gamma \mathcal{K}_{2}+\gamma (\gamma -1)\mathcal{K}%
_{1}+\mathcal{K}_{0}}{\gamma \left( \gamma -\beta \right) \left( \gamma
-\alpha \right) }.
\end{eqnarray*}
\end{theorem}

Proof. We prove the theorem only for $n\geq 0$ because of Theorem \ref%
{theorem:fvgxdfsxczsaer}. We prove (\ref{equation:mousxzadfgytc}). By the
assumption, the characteristic equation of (\ref{equation:yusoabpsmnbscv})
is $x^{3}-x^{2}-x-1=0$ and the roots of it are $\alpha ,\beta $ and $\gamma
. $ So it's general solution is given by%
\begin{equation*}
\mathcal{T}_{n}=A_{1}\alpha ^{n}+B_{1}\beta ^{n}+C_{1}\gamma ^{n}.
\end{equation*}%
Using initial condition which is given in Definition \ref%
{definition:tyfsdczxawuqp}, and also applying lineer algebra operations, we
obtain the matrices $A_{1},B_{1},C_{1}$ as desired. This gives the formula
for $\mathcal{T}_{n}.$

Similarly we have the formula (\ref{equat:fvbxdszuqwsazx}). 
\endproof%

The well known Binet formulas for Tribonacci and Tribonacci-Lucas numbers
are given in (\ref{equat:mnopcvbedcxzsa})\ and (\ref{equation:cfrdcsxszouea}%
) respectively. But, we will obtain these functions in terms of Tribonacci
and Tribonacci-Lucas matrix sequences as a consequence of Theorems \ref%
{theorem:hbnvxgfdtreopas} and \ref{theorem:sacasdmnbhgvc}. To do this, we
will give the formulas for these numbers by means of the related matrix
sequences. In fact, in the proof of next corollary, we will just compare the
linear combination of the 2nd row and 1st column entries of the matrices.

\begin{corollary}
\label{corollary:fgvbxdzfserquevbdf}For every integers $n,$ the Binet's
formulas for Tribonacci and Tribonacci-Lucas numbers are given as%
\begin{eqnarray*}
T_{n} &=&\frac{\alpha ^{n+1}}{\left( \alpha -\gamma \right) \left( \alpha
-\beta \right) }+\frac{\beta ^{n+1}}{\left( \beta -\gamma \right) \left(
\beta -\alpha \right) }+\frac{\gamma ^{n+1}}{\left( \gamma -\beta \right)
\left( \gamma -\alpha \right) }, \\
K_{n} &=&\alpha ^{n}+\beta ^{n}+\gamma ^{n}.
\end{eqnarray*}
\end{corollary}

\textit{Proof.} From Theorem \ref{theorem:sacasdmnbhgvc}, we have 
\begin{eqnarray*}
\mathcal{T}_{n} &=&A_{1}\alpha ^{n}+B_{1}\beta ^{n}+C_{1}\gamma ^{n} \\
&=&\frac{\alpha \mathcal{T}_{2}+\alpha (\alpha -1)\mathcal{T}_{1}+\mathcal{T}%
_{0}}{\alpha \left( \alpha -\gamma \right) \left( \alpha -\beta \right) }%
\alpha ^{n}+\frac{\beta \mathcal{T}_{2}+\beta (\beta -1)\mathcal{T}_{1}+%
\mathcal{T}_{0}}{\beta \left( \beta -\gamma \right) \left( \beta -\alpha
\right) }\beta ^{n} \\
&&+\frac{\gamma \mathcal{T}_{2}+\gamma (\gamma -1)\mathcal{T}_{1}+\mathcal{T}%
_{0}}{\gamma \left( \gamma -\beta \right) \left( \gamma -\alpha \right) }%
\gamma ^{n} \\
&=&\frac{\alpha ^{n-1}}{\left( \alpha -\gamma \right) \left( \alpha -\beta
\right) }\left( 
\begin{array}{ccc}
\alpha ^{3} & \alpha \left( \alpha +1\right) & \alpha ^{2} \\ 
\alpha ^{2} & \alpha +1 & \alpha \\ 
\alpha & \alpha \left( \alpha -1\right) & 1%
\end{array}%
\right) +\frac{\beta ^{n-1}}{\left( \beta -\gamma \right) \left( \beta
-\alpha \right) }\left( 
\begin{array}{ccc}
\beta ^{3} & \beta \left( \beta +1\right) & \beta ^{2} \\ 
\beta ^{2} & \beta +1 & \beta \\ 
\beta & \beta \left( \beta -1\right) & 1%
\end{array}%
\right) \\
&&+\frac{\gamma ^{n-1}}{\left( \gamma -\beta \right) \left( \gamma -\alpha
\right) }\left( 
\begin{array}{ccc}
\gamma ^{3} & \gamma \left( \gamma +1\right) & \gamma ^{2} \\ 
\gamma ^{2} & \gamma +1 & \gamma \\ 
\gamma & \gamma \left( \gamma -1\right) & 1%
\end{array}%
\right)
\end{eqnarray*}%
By Theorem \ref{theorem:hbnvxgfdtreopas}, we know that 
\begin{equation*}
\mathcal{T}_{n}=\left( 
\begin{array}{ccc}
T_{n+1} & T_{n}+T_{n-1} & T_{n} \\ 
T_{n} & T_{n-1}+T_{n-2} & T_{n-1} \\ 
T_{n-1} & T_{n-2}+T_{n-3} & T_{n-2}%
\end{array}%
\right) .
\end{equation*}%
Now, if we compare the 2nd row and 1st column entries with the matrices in
the above two equations, then we obtain%
\begin{eqnarray*}
T_{n} &=&\frac{\alpha ^{n-1}\alpha ^{2}}{\left( \alpha -\gamma \right)
\left( \alpha -\beta \right) }+\frac{\beta ^{n-1}\beta ^{2}}{\left( \beta
-\gamma \right) \left( \beta -\alpha \right) }+\frac{\gamma ^{n-1}\gamma ^{2}%
}{\left( \gamma -\beta \right) \left( \gamma -\alpha \right) } \\
&=&\frac{\alpha ^{n+1}}{\left( \alpha -\gamma \right) \left( \alpha -\beta
\right) }+\frac{\beta ^{n+1}}{\left( \beta -\gamma \right) \left( \beta
-\alpha \right) }+\frac{\gamma ^{n+1}}{\left( \gamma -\beta \right) \left(
\gamma -\alpha \right) }.
\end{eqnarray*}%
From Therem \ref{theorem:sacasdmnbhgvc}, we obtain%
\begin{eqnarray*}
\mathcal{K}_{n} &=&A_{2}\alpha ^{n}+B_{2}\beta ^{n}+C_{2}\gamma ^{n} \\
&=&\frac{\alpha \mathcal{K}_{2}+\alpha (\alpha -1)\mathcal{K}_{1}+\mathcal{K}%
_{0}}{\alpha \left( \alpha -\gamma \right) \left( \alpha -\beta \right) }%
\alpha ^{n}+\frac{\beta \mathcal{K}_{2}+\beta (\beta -1)\mathcal{K}_{1}+%
\mathcal{K}_{0}}{\beta \left( \beta -\gamma \right) \left( \beta -\alpha
\right) }\beta ^{n} \\
&&+\frac{\gamma \mathcal{K}_{2}+\gamma (\gamma -1)\mathcal{K}_{1}+\mathcal{K}%
_{0}}{\gamma \left( \gamma -\beta \right) \left( \gamma -\alpha \right) }%
\gamma ^{n} \\
&=&\frac{\alpha ^{n-1}}{\left( \alpha -\gamma \right) \left( \alpha -\beta
\right) }\left( 
\begin{array}{ccc}
3\alpha ^{2}+4\alpha +1 & 4\alpha ^{2}+2 & \alpha ^{2}+2\alpha +3 \\ 
\alpha ^{2}+2\alpha +3 & 2\alpha ^{2}+2\alpha -2 & 3\alpha ^{2}-2\alpha -1
\\ 
3\alpha ^{2}-2\alpha -1 & -2\alpha ^{2}+4\alpha +4 & -\alpha ^{2}+4\alpha -1%
\end{array}%
\right) \\
&&+\frac{\beta ^{n-1}}{\left( \beta -\gamma \right) \left( \beta -\alpha
\right) }\left( 
\begin{array}{ccc}
3\beta ^{2}+4\beta +1 & 4\beta ^{2}+2 & \beta ^{2}+2\beta +3 \\ 
\beta ^{2}+2\beta +3 & 2\beta ^{2}+2\beta -2 & 3\beta ^{2}-2\beta -1 \\ 
3\beta ^{2}-2\beta -1 & -2\beta ^{2}+4\beta +4 & -\beta ^{2}+4\beta -1%
\end{array}%
\right) \\
&&+\frac{\gamma ^{n-1}}{\left( \gamma -\beta \right) \left( \gamma -\alpha
\right) }\left( 
\begin{array}{ccc}
3\gamma ^{2}+4\gamma +1 & 4\gamma ^{2}+2 & \gamma ^{2}+2\gamma +3 \\ 
\gamma ^{2}+2\gamma +3 & 2\gamma ^{2}+2\gamma -2 & 3\gamma ^{2}-2\gamma -1
\\ 
3\gamma ^{2}-2\gamma -1 & -2\gamma ^{2}+4\gamma +4 & -\gamma ^{2}+4\gamma -1%
\end{array}%
\right) .
\end{eqnarray*}%
By Theorem \ref{theorem:hbnvxgfdtreopas}, we know that 
\begin{equation*}
\mathcal{K}_{n}=\left( 
\begin{array}{ccc}
K_{n+1} & K_{n}+K_{n-1} & K_{n} \\ 
K_{n} & K_{n-1}+K_{n-2} & K_{n-1} \\ 
K_{n-1} & K_{n-2}+K_{n-3} & K_{n-2}%
\end{array}%
\right) .
\end{equation*}%
Now, if we compare the 2nd row and 1st column entries with the matrices in
the above last two equations, then we obtain%
\begin{equation*}
K_{n}=\frac{\alpha ^{n-1}(\alpha ^{2}+2\alpha +3)}{\left( \alpha -\gamma
\right) \left( \alpha -\beta \right) }+\frac{\beta ^{n-1}(\beta ^{2}+2\beta
+3)}{\left( \beta -\gamma \right) \left( \beta -\alpha \right) }+\frac{%
\gamma ^{n-1}(\gamma ^{2}+2\gamma +3)}{\left( \gamma -\beta \right) \left(
\gamma -\alpha \right) }.
\end{equation*}%
Using the relations, $\alpha +\beta +\gamma =1,$ $\alpha \beta \gamma =1$
and considering $\alpha ,\beta $ and $\gamma $ are the roots the equation $%
x^{3}-x^{2}-x-1=0,$ we obtain 
\begin{eqnarray*}
\frac{\alpha ^{2}+2\alpha +3}{\left( \alpha -\gamma \right) \left( \alpha
-\beta \right) } &=&\frac{\alpha ^{2}+2\alpha +3}{\alpha ^{2}-\alpha \beta
-\alpha \gamma +\beta \gamma }=\frac{\alpha }{\alpha }\frac{\alpha
^{2}+2\alpha +3}{\alpha ^{2}+\alpha (-\beta -\gamma )+\beta \gamma } \\
&=&\frac{\alpha }{\alpha }\frac{(\alpha ^{2}+2\alpha +3)}{\alpha ^{2}+\alpha
(-\beta -\gamma )+\beta \gamma }=\frac{(\alpha ^{2}+2\alpha +3)\alpha }{%
\alpha ^{3}+\alpha ^{2}(\alpha -1)+1} \\
&=&\frac{(\alpha ^{2}+2\alpha +3)\alpha }{2\alpha ^{3}-\alpha ^{2}+1}=\frac{%
(\alpha ^{2}+2\alpha +3)\alpha }{2(\alpha ^{2}+\alpha +1)-\alpha ^{2}+1} \\
&=&\frac{(\alpha ^{2}+2\alpha +3)\alpha }{(\alpha ^{2}+2\alpha +3)}=\alpha ,
\\
\frac{\beta ^{2}+2\beta +3}{\left( \beta -\gamma \right) \left( \beta
-\alpha \right) } &=&\frac{\beta ^{2}+2\beta +3}{\beta ^{2}-\alpha \beta
+\alpha \gamma -\beta \gamma }=\beta , \\
\frac{\gamma ^{2}+2\gamma +3}{\left( \gamma -\beta \right) \left( \gamma
-\alpha \right) } &=&\frac{\gamma ^{2}+2\gamma +3}{\gamma ^{2}+\alpha \beta
-\alpha \gamma -\beta \gamma }=\gamma .
\end{eqnarray*}%
So finally we conclude that 
\begin{equation*}
K_{n}=\alpha ^{n}+\beta ^{n}+\gamma ^{n}
\end{equation*}%
as required. 
\endproof%

Now, we present summation formulas for Tribonacci and Tribonacci-Lucas
matrix sequences.

\begin{theorem}
For $m>j\geq 0,$ we have%
\begin{equation}
\sum_{i=0}^{n-1}\mathcal{T}_{mi+j}=\frac{\mathcal{T}_{mn+m+j}+\mathcal{T}%
_{mn-m+j}+(1-K_{m})\mathcal{T}_{mn+j}}{K_{m}-K_{-m}}-\frac{\mathcal{T}_{m+j}+%
\mathcal{T}_{j-m}+(1-K_{m})\mathcal{T}_{j}}{K_{m}-K_{-m}}
\label{equati:bnosdczxdaer}
\end{equation}
and%
\begin{equation}
\sum_{i=0}^{n-1}\mathcal{K}_{mi+j}=\frac{\mathcal{K}_{mn+m+j}+\mathcal{K}%
_{mn-m+j}+(1-K_{m})\mathcal{K}_{mn+j}}{K_{m}-K_{-m}}-\frac{\mathcal{K}_{m+j}+%
\mathcal{K}_{j-m}+(1-K_{m})\mathcal{K}_{j}}{K_{m}-K_{-m}}.
\label{equat:bnyugfdsxcza}
\end{equation}
\end{theorem}

\textit{Proof. Note that} 
\begin{eqnarray*}
\sum_{i=0}^{n-1}\mathcal{T}_{mi+j} &=&\sum_{i=0}^{n-1}(A_{1}\alpha
^{mi+j}+B_{1}\beta ^{mi+j}+C_{1}\gamma ^{mi+j}) \\
&=&A_{1}\alpha ^{j}\left( \frac{\alpha ^{mn}-1}{\alpha ^{m}-1}\right)
+B_{1}\beta ^{j}\left( \frac{\beta ^{mn}-1}{\beta ^{m}-1}\right)
+C_{1}\gamma ^{j}\left( \frac{\gamma ^{mn}-1}{\gamma ^{m}-1}\right)
\end{eqnarray*}%
and%
\begin{eqnarray*}
\sum_{i=0}^{n-1}\mathcal{K}_{mi+j} &=&\sum_{i=0}^{n-1}(A_{2}\alpha
^{mi+j}+B_{2}\beta ^{mi+j}+C_{2}\gamma ^{mi+j}) \\
&=&A_{2}\alpha ^{j}\left( \frac{\alpha ^{mn}-1}{\alpha ^{m}-1}\right)
+B_{2}\beta ^{j}\left( \frac{\beta ^{mn}-1}{\beta ^{m}-1}\right)
+C_{2}\gamma ^{j}\left( \frac{\gamma ^{mn}-1}{\gamma ^{m}-1}\right) .
\end{eqnarray*}%
Simplifying and rearranging the last equalities in the last two expression
imply (\ref{equati:bnosdczxdaer}) and (\ref{equat:bnyugfdsxcza}) as
required. 
\endproof%

As in Corollary \ref{corollary:fgvbxdzfserquevbdf}, in the proof of next
Corollary, we just compare the linear combination of the 2nd row and 1st
column entries of the relevant matrices.

\begin{corollary}
For $m>j\geq 0,$ we have%
\begin{equation}
\sum_{i=0}^{n-1}T_{mi+j}=\frac{T_{mn+m+j}+T_{mn-m+j}+(1-K_{m})T_{mn+j}}{%
K_{m}-K_{-m}}-\frac{T_{m+j}+T_{j-m}+(1-K_{m})T_{j}}{K_{m}-K_{-m}}
\end{equation}%
and%
\begin{equation}
\sum\limits_{i=0}^{n-1}K_{mi+j}=\frac{K_{mn+m+j}+K_{mn-m+j}+(1-K_{m})K_{mn+j}%
}{K_{m}-K_{-m}}-\frac{K_{m+j}+K_{j-m}+(1-K_{m})K_{j}}{K_{m}-K_{-m}}
\end{equation}
\end{corollary}

Note that using the above Corollary we obtain the following well known
formulas (taking $m=1,j=0$):%
\begin{equation*}
\sum_{i=0}^{n-1}T_{i}=\frac{T_{n+2}-T_{n}-1}{2}\text{ \ and }%
\sum_{i=0}^{n-1}K_{i}=\frac{K_{n+2}-K_{n}}{2}.\text{ \ }
\end{equation*}

We now give generating functions of $\mathcal{T}$ and $\mathcal{K}$.

\begin{theorem}
\label{theorem:avsonmbvxcfd}The generating function for the Tribonacci and
Tribonacci-Lucas matrix sequences are given as%
\begin{equation*}
\sum_{n=0}^{\infty }\mathcal{T}_{n}x^{n}=\frac{1}{1-x-x^{2}-x^{3}}\left( 
\begin{array}{ccc}
1 & x+x^{2} & x \\ 
x & 1-x & x^{2} \\ 
x^{2} & x-x^{2} & 1-x-x^{2}%
\end{array}%
\right)
\end{equation*}%
and%
\begin{equation*}
\sum_{n=0}^{\infty }\mathcal{K}_{n}x^{n}=\frac{1}{1-x-x^{2}-x^{3}}\left( 
\begin{array}{ccc}
1+2x+3x^{2} & 2+2x-2x^{2} & 3-2x-x^{2} \\ 
3-2x-x^{2} & -2+4x+4x^{2} & -1+4x-x^{2} \\ 
-1+4x-x^{2} & 4-6x & -1+5x^{2}%
\end{array}%
\right)
\end{equation*}%
respectively.
\end{theorem}

\textit{Proof. }We prove the Tribonacci case. Suppose that $%
g(x)=\sum_{n=0}^{\infty }\mathcal{T}_{n}x^{n}$ is the generating function
for the sequence $\{\mathcal{T}_{n}\}_{n\geq 0}.$ Then, using Definition \ref%
{definition:tyfsdczxawuqp}, we obtain%
\begin{eqnarray*}
g(x) &=&\sum_{n=0}^{\infty }\mathcal{T}_{n}x^{n}=\mathcal{T}_{0}+\mathcal{T}%
_{1}x+\mathcal{T}_{2}x^{2}+\sum_{n=3}^{\infty }\mathcal{T}_{n}x^{n} \\
&=&\mathcal{T}_{0}+\mathcal{T}_{1}x+\mathcal{T}_{2}x^{2}+\sum_{n=3}^{\infty
}(\mathcal{T}_{n-1}+\mathcal{T}_{n-2}+\mathcal{T}_{n-3})x^{n} \\
&=&\mathcal{T}_{0}+\mathcal{T}_{1}x+\mathcal{T}_{2}x^{2}+\sum_{n=3}^{\infty }%
\mathcal{T}_{n-1}x^{n}+\sum_{n=3}^{\infty }\mathcal{T}_{n-2}x^{n}+%
\sum_{n=3}^{\infty }\mathcal{T}_{n-3}x^{n} \\
&=&\mathcal{T}_{0}+\mathcal{T}_{1}x+\mathcal{T}_{2}x^{2}-\mathcal{T}_{0}x-%
\mathcal{T}_{1}x^{2}-\mathcal{T}_{0}x^{2}+x\sum_{n=0}^{\infty }\mathcal{T}%
_{n}x^{n}+x^{2}\sum_{n=0}^{\infty }\mathcal{T}_{n}x^{n}+x^{3}\sum_{n=0}^{%
\infty }\mathcal{T}_{n}x^{n} \\
&=&\mathcal{T}_{0}+\mathcal{T}_{1}x+\mathcal{T}_{2}x^{2}-\mathcal{T}_{0}x-%
\mathcal{T}_{1}x^{2}-\mathcal{T}_{0}x^{2}+xg(x)+x^{2}g(x)+x^{3}g(x).
\end{eqnarray*}%
Rearranging above equation, we get%
\begin{equation*}
g(x)=\frac{\mathcal{T}_{0}+(\mathcal{T}_{1}-\mathcal{T}_{0})x+(\mathcal{T}%
_{2}-\mathcal{T}_{1}-\mathcal{T}_{0})x^{2}}{1-x-x^{2}-x^{3}}.
\end{equation*}%
which equals the $\sum_{n=0}^{\infty }\mathcal{T}_{n}x^{n}$ in the Theorem.
This completes the proof.

Tribonacci-Lucas case can be proved similarly. 
\endproof%

The well known generating functions for Tribonacci and Tribonacci-Lucas
numbers are as in (\ref{equation:yugdfvxbgsopqac}). However, we will obtain
these functions in terms of Tribonacci and Tribonacci-Lucas matrix sequences
as a consequence of Theorem \ref{theorem:avsonmbvxcfd}. To do this, we will
again compare the the 2nd row and 1st column entries with the matrices in
Theorem \ref{theorem:avsonmbvxcfd}. Thus we have the following corollary.

\begin{corollary}
\bigskip The generating functions for the Tribonacci sequence $%
\{T_{n}\}_{n\geq 0}$ and Tribonacci-Lucas sequence $\{K_{n}\}_{n\geq 0}$ are
given as%
\begin{equation*}
\sum_{n=0}^{\infty }T_{n}x^{n}=\frac{x}{1-x-x^{2}-x^{3}}\text{\ and \ }%
\sum_{n=0}^{\infty }K_{n}x^{n}=\frac{3-2x-x^{2}}{1-x-x^{2}-x^{3}}.
\end{equation*}%
respectively.
\end{corollary}

\section{Relation Between Tribonacci and Tribonacci-Lucas Matrix Sequences}

The following theorem shows that there always exist interrelation between
Tribonacci and Tribonacci-Lucas matrix sequences.

\begin{theorem}
\label{theorem:oxcvztysufrewdxs}For the matrix sequences $\{\mathcal{T}%
_{n}\} $ and $\{\mathcal{K}_{n}\},\ $we have the following identities.

\begin{description}
\item[(a)] $\mathcal{K}_{n}=3\mathcal{T}_{n+1}-2\mathcal{T}_{n}-\mathcal{T}%
_{n-1},$

\item[(b)] $\mathcal{K}_{n}=\mathcal{T}_{n}+2\mathcal{T}_{n-1}+3\mathcal{T}%
_{n-2},$

\item[(c)] $\mathcal{K}_{n}=4\mathcal{T}_{n+1}-\mathcal{T}_{n}-\mathcal{T}%
_{n+2},$

\item[(d)] $\mathcal{K}_{n}=-\mathcal{T}_{n+2}+4\mathcal{T}_{n+1}-\mathcal{T}%
_{n},$

\item[(e)] $\mathcal{T}_{n}=\frac{1}{22}(5\mathcal{K}_{n+2}-3\mathcal{K}%
_{n+1}-4\mathcal{K}_{n})$
\end{description}
\end{theorem}

\textit{Proof.} From (\ref{equation:gfdvbxczvsadou}), (\ref%
{equation:fvcdxsartewqa}) and (\ref{equation:yuhtgsdafcscvb}), (a), (b) and
(c) follow. It is easy to show that $K_{n}=-T_{n+2}+4T_{n+1}-T_{n}$ and $%
22T_{n}=5K_{n+2}-3K_{n+1}-4K_{n}$ using Binet formulas of the numbers $T_{n}$
and $K_{n},$ so now (d) and (e) follow. 
\endproof%

\begin{lemma}
\label{lemma:uygbcxfdzreawqops}For all non-negative integers $m$ and $n,\ $%
we have the following identities.

\begin{description}
\item[(a)] $\mathcal{K}_{0}\mathcal{T}_{n}=\mathcal{T}_{n}\mathcal{K}_{0}=%
\mathcal{K}_{n},$

\item[(b)] $\mathcal{T}_{0}\mathcal{K}_{n}=\mathcal{K}_{n}\mathcal{T}_{0}=%
\mathcal{K}_{n}.$
\end{description}
\end{lemma}

\textit{Proof.} Identities can be established easily. Note that to show (a)
we need to use all the relations (\ref{equation:gfdvbxczvsadou}), (\ref%
{equation:fvcdxsartewqa}) and (\ref{equation:yuhtgsdafcscvb}). 
\endproof%

Next Corollary gives another relation between the numbers $T_{n}$ and $K_{n}$
and also the matrices $\mathcal{T}_{n}$ and $\mathcal{K}_{n}$.

\begin{corollary}
We have the following identities.

\begin{description}
\item[(a)] $T_{n}=\frac{1}{22}(K_{n}+5K_{n-1}+2K_{n+1}),$

\item[(b)] $\mathcal{T}_{n}=\frac{1}{22}(\mathcal{K}_{n}+5\mathcal{K}_{n-1}+2%
\mathcal{K}_{n+1}).$
\end{description}
\end{corollary}

\textit{Proof.} From Lemma \ref{lemma:uygbcxfdzreawqops} (a), we know that $%
\mathcal{K}_{0}\mathcal{T}_{n}=\mathcal{K}_{n}.$ To show (a), use Theorem %
\ref{theorem:hbnvxgfdtreopas} for the matrix $\mathcal{T}_{n}$ and calculate
the matrix operation $\mathcal{K}_{0}^{-1}\mathcal{K}_{n}$ and then compare
the 2nd row and 1st column entries with the matrices $\mathcal{T}_{n}$ and $%
\mathcal{K}_{0}^{-1}\mathcal{K}_{n}.$ Now (b) follows from (a). 
\endproof%

To prove the following Theorem we need the next Lemma.

\begin{lemma}
\label{lemma:gbnhytrdvopxdser}Let $A_{1},B_{1},C_{1};A_{2},B_{2},C_{2}$ as
in Theorem \ref{theorem:sacasdmnbhgvc}. Then the following relations hold:%
\begin{eqnarray*}
A_{1}^{2} &=&A_{1},\text{ }B_{1}^{2}=B_{1},\text{ }C_{1}^{2}=C_{1}, \\
A_{1}B_{1} &=&B_{1}A_{1}=A_{1}C_{1}=C_{1}A_{1}=C_{1}B_{1}=B_{1}C_{1}=\left(
0\right) , \\
A_{2}B_{2} &=&B_{2}A_{2}=A_{2}C_{2}=C_{2}A_{2}=C_{2}B_{2}=B_{2}C_{2}=\left(
0\right) .
\end{eqnarray*}
\end{lemma}

\textit{Proof. }Using $\alpha +\beta +\gamma =1,$ $\alpha \beta +\alpha
\gamma +\beta \gamma =-1$ and $\alpha \beta \gamma =1,$ required equalities
can be established by matrix calculations. 
\endproof%

\begin{theorem}
\label{theorem:fvcxoueaxsdwqa}For all non-negative integers $m$ and $n,\ $we
have the following identities.

\begin{description}
\item[(a)] $\mathcal{T}_{m}\mathcal{T}_{n}=\mathcal{T}_{m+n}=\mathcal{T}_{n}%
\mathcal{T}_{m},$

\item[(b)] $\mathcal{T}_{m}\mathcal{K}_{n}=\mathcal{K}_{n}\mathcal{T}_{m}=%
\mathcal{K}_{m+n},$

\item[(c)] $\mathcal{K}_{m}\mathcal{K}_{n}=\mathcal{K}_{n}\mathcal{K}_{m}=9%
\mathcal{T}_{m+n+2}-12\mathcal{T}_{m+n+1}-2\mathcal{T}_{m+n}+4\mathcal{T}%
_{m+n-1}+\mathcal{T}_{m+n-2},$

\item[(d)] $\mathcal{K}_{m}\mathcal{K}_{n}=\mathcal{K}_{n}\mathcal{K}_{m}=%
\mathcal{T}_{m+n}+4\mathcal{T}_{m+n-1}+10\mathcal{T}_{m+n-2}+12\mathcal{T}%
_{m+n-3}+\allowbreak 9\mathcal{T}_{m+n-4},$

\item[(e)] $\mathcal{K}_{m}\mathcal{K}_{n}=\mathcal{K}_{n}\mathcal{K}_{m}=%
\mathcal{T}_{m+n}-8\mathcal{T}_{m+n+1}+18\mathcal{T}_{m+n+2}-8\mathcal{T}%
_{m+n+3}+\mathcal{T}_{m+n+4}.$
\end{description}
\end{theorem}

\textit{Proof.} \ 

\begin{description}
\item[(a)] Using Lemma \ref{lemma:gbnhytrdvopxdser} we obtain 
\begin{eqnarray*}
\mathcal{T}_{m}\mathcal{T}_{n} &=&(A_{1}\alpha ^{m}+B_{1}\beta
^{m}+C_{1}\gamma ^{m})(A_{1}\alpha ^{n}+B_{1}\beta ^{n}+C_{1}\gamma ^{n}) \\
&=&A_{1}^{2}\alpha ^{m+n}+B_{1}^{2}\beta ^{m+n}+C_{1}^{2}\gamma
^{m+n}+A_{1}\allowbreak B_{1}\alpha ^{m}\beta ^{n}+B_{1}A_{1}\alpha
^{n}\beta ^{m} \\
&&+A_{1}C_{1}\alpha ^{m}\gamma ^{n}+C_{1}A_{1}\alpha ^{n}\gamma
^{m}+B_{1}C_{1}\beta ^{m}\allowbreak \gamma ^{n}+C_{1}B_{1}\beta ^{n}\gamma
^{m} \\
&=&A_{1}\alpha ^{m+n}+B_{1}\beta ^{m+n}+C_{1}\gamma ^{m+n} \\
&=&\mathcal{T}_{m+n}.
\end{eqnarray*}

\item[(b)] By Lemma \ref{lemma:uygbcxfdzreawqops}, we have%
\begin{equation*}
\mathcal{T}_{m}\mathcal{K}_{n}=\mathcal{T}_{m}\mathcal{T}_{n}\mathcal{K}_{0}.
\end{equation*}%
Now from (a) and again by Lemma \ref{lemma:uygbcxfdzreawqops} we obtain $%
\mathcal{T}_{m}\mathcal{K}_{n}=\mathcal{T}_{m+n}\mathcal{K}_{0}=\mathcal{K}%
_{m+n}$.

It can be shown similarly that $\mathcal{K}_{n}\mathcal{T}_{m}=\mathcal{K}%
_{m+n}.$

\item[(c)] Using (a) and Theorem \ref{theorem:oxcvztysufrewdxs} (a) we
obtain 
\begin{eqnarray*}
\mathcal{K}_{m}\mathcal{K}_{n} &=&(3\mathcal{T}_{m+1}-2\mathcal{T}_{m}-%
\mathcal{T}_{m-1})(3\mathcal{T}_{n+1}-2\mathcal{T}_{n}-\mathcal{T}_{n-1}) \\
&=&2\mathcal{T}_{n}\mathcal{T}_{m-1}-6\mathcal{T}_{n}\mathcal{T}_{m+1}+2%
\mathcal{T}_{m}\mathcal{T}_{n-1}-6\mathcal{T}_{m}\mathcal{T}_{n+1} \\
&&+\allowbreak 4\mathcal{T}_{m}\mathcal{T}_{n}+\mathcal{T}_{m-1}\mathcal{T}%
_{n-1}-3\mathcal{T}_{m-1}\mathcal{T}_{n+1}-3\mathcal{T}_{m+1}\mathcal{T}%
_{n-1}+\allowbreak 9\mathcal{T}_{m+1}\mathcal{T}_{n+1} \\
&=&2\mathcal{T}_{m+n-1}-6\mathcal{T}_{m+n+1}+2\mathcal{T}_{m+n-1}-6\mathcal{T%
}_{m+n+1}+\allowbreak 4\mathcal{T}_{m+n}+\mathcal{T}_{m+n-2}-3\mathcal{T}%
_{m+n} \\
&&-3\mathcal{T}_{m+n}+\allowbreak 9\mathcal{T}_{m+n+2} \\
&=&9\mathcal{T}_{m+n+2}-12\mathcal{T}_{m+n+1}-2\mathcal{T}_{m+n}+4\mathcal{T}%
_{m+n-1}+\mathcal{T}_{m+n-2}
\end{eqnarray*}%
$\allowbreak $

It can be shown similarly that $\mathcal{K}_{n}\mathcal{K}_{m}=9\mathcal{T}%
_{m+n+2}-12\mathcal{T}_{m+n+1}-2\mathcal{T}_{m+n}+4\mathcal{T}_{m+n-1}+%
\mathcal{T}_{m+n-2}.$

The remaining of identities can be proved by considering again (a) and
Theorem \ref{theorem:oxcvztysufrewdxs}. 
\endproof%
\end{description}

Comparing matrix entries and using Teorem \ref{theorem:hbnvxgfdtreopas} we
have next result.

\begin{corollary}
For Tribonacci and Tribonacci-Lucas numbers, we have the following
identities:

\begin{description}
\item[(a)] $T_{m+n}=T_{m}T_{n+1}+T_{n}\left( T_{m-1}+T_{m-2}\right)
+T_{m-1}T_{n-1}$

\item[(b)] $K_{m+n}=T_{m}K_{n+1}+K_{n}\left( T_{m-1}+T_{m-2}\right)
+K_{n-1}T_{m-1}$

\item[(c)] $K_{m}K_{n+1}+K_{n}\left( K_{m-1}+K_{m-2}\right)
+K_{m-1}K_{n-1}=9T_{m+n+2}-12T_{m+n+1}-2T_{m+n}+4T_{m+n-1}+T_{m+n-2}$

\item[(d)] $K_{m}K_{n+1}+K_{n}\left( K_{m-1}+K_{m-2}\right)
+K_{m-1}K_{n-1}=\allowbreak
T_{m+n}+4T_{m+n-1}+10T_{m+n-2}+12T_{m+n-3}+9T_{m+n-4}$

\item[(e)] $K_{m}K_{n+1}+K_{n}\left( K_{m-1}+K_{m-2}\right)
+K_{m-1}K_{n-1}=T_{m+n}-8T_{m+n+1}+18T_{m+n+2}-8T_{m+n+3}+T_{m+n+4}$
\end{description}
\end{corollary}

\textit{Proof.}

\begin{description}
\item[(a)] From Theorem \ref{theorem:fvcxoueaxsdwqa} we know that $\mathcal{T%
}_{m}\mathcal{T}_{n}=\mathcal{T}_{m+n}.$ Using Theorem \ref%
{theorem:hbnvxgfdtreopas}, we can write this result as%
\begin{eqnarray*}
&&\left( 
\begin{array}{ccc}
T_{m+1} & T_{m}+T_{m-1} & T_{m} \\ 
T_{m} & T_{m-1}+T_{m-2} & T_{m-1} \\ 
T_{m-1} & T_{m-2}+T_{m-3} & T_{m-2}%
\end{array}%
\right) \left( 
\begin{array}{ccc}
T_{n+1} & T_{n}+T_{n-1} & T_{n} \\ 
T_{n} & T_{n-1}+T_{n-2} & T_{n-1} \\ 
T_{n-1} & T_{n-2}+T_{n-3} & T_{n-2}%
\end{array}%
\right) \\
&=&\left( 
\begin{array}{ccc}
T_{m+n+1} & T_{m+n}+T_{m+n-1} & T_{m+n} \\ 
T_{m+n} & T_{m+n-1}+T_{m+n-2} & T_{m+n-1} \\ 
T_{m+n-1} & T_{m+n-2}+T_{m+n-3} & T_{m+n-2}%
\end{array}%
\right) .
\end{eqnarray*}%
Now, by multiplying the left-side matrices and then by comparing the 2nd
rows and 1st columns entries, we get the required identity in (a).

The remaining of identities can be proved by considering again Theorems \ref%
{theorem:fvcxoueaxsdwqa}\ and \ref{theorem:hbnvxgfdtreopas}. 
\endproof%
\end{description}

The next two theorems provide us the convenience to obtain the powers of
Tribonacci and Tribonacci-Lucas matrix sequences.

\begin{theorem}
\label{theorem:cvxhyugdswazxsa}For non-negatif integers $m,n$ and $r$ with $%
n\geq r,$ the following identities hold:

\begin{description}
\item[(a)] $\mathcal{T}_{n}^{m}=\mathcal{T}_{mn},$

\item[(b)] $\mathcal{T}_{n+1}^{m}=\mathcal{T}_{1}^{m}\mathcal{T}_{mn},$

\item[(c)] $\mathcal{T}_{n-r}\mathcal{T}_{n+r}=\mathcal{T}_{n}^{2}=\mathcal{T%
}_{2}^{n}.$
\end{description}
\end{theorem}

\textit{Proof.} \ \ 

\begin{description}
\item[(a)] We can write $\mathcal{T}_{n}^{m}$ as%
\begin{equation*}
\mathcal{T}_{n}^{m}=\mathcal{T}_{n}\mathcal{T}_{n}...\mathcal{T}_{n}\text{ (}%
m\text{ times).}
\end{equation*}%
Using Theorem \ref{theorem:fvcxoueaxsdwqa} (a) iteratively, we obtain the
required result:%
\begin{eqnarray*}
\mathcal{T}_{n}^{m} &=&\underset{m\text{ times}}{\underbrace{\mathcal{T}_{n}%
\mathcal{T}_{n}...\mathcal{T}_{n}}} \\
&=&\mathcal{T}_{2n}\underset{m-1\text{ times}}{\underbrace{\mathcal{T}_{n}%
\mathcal{T}_{n}...\mathcal{T}_{n}}} \\
&=&\mathcal{T}_{3n}\underset{m-2\text{ times}}{\underbrace{\mathcal{T}_{n}%
\mathcal{T}_{n}...\mathcal{T}_{n}}} \\
&&\vdots \\
&=&\mathcal{T}_{(m-1)n}\mathcal{T}_{n} \\
&=&\mathcal{T}_{mn}.
\end{eqnarray*}

\item[(b)] As a similar approach in (a) we have%
\begin{equation*}
\mathcal{T}_{n+1}^{m}=\mathcal{T}_{n+1}.\mathcal{T}_{n+1}...\mathcal{T}%
_{n+1}=\mathcal{T}_{m(n+1)}=\mathcal{T}_{m}\mathcal{T}_{mn}=\mathcal{T}_{1}%
\mathcal{T}_{m-1}\mathcal{T}_{mn}.
\end{equation*}%
Using Theorem \ref{theorem:fvcxoueaxsdwqa} (a), we can write iteratively $%
\mathcal{T}_{m}=\mathcal{T}_{1}\mathcal{T}_{m-1},$ $\mathcal{T}_{m-1}=%
\mathcal{T}_{1}\mathcal{T}_{m-2},$ ..., $\mathcal{T}_{2}=\mathcal{T}_{1}%
\mathcal{T}_{1}.$ Now it follows that 
\begin{equation*}
\mathcal{T}_{n+1}^{m}=\underset{m\text{ times}}{\underbrace{\mathcal{T}_{1}%
\mathcal{T}_{1}...\mathcal{T}_{1}}}\mathcal{T}_{mn}=\mathcal{T}_{1}^{m}%
\mathcal{T}_{mn}.
\end{equation*}

\item[(c)] Theorem \ref{theorem:fvcxoueaxsdwqa} (a) gives%
\begin{equation*}
\mathcal{T}_{n-r}\mathcal{T}_{n+r}=\mathcal{T}_{2n}=\mathcal{T}_{n}\mathcal{T%
}_{n}=\mathcal{T}_{n}^{2}
\end{equation*}%
and also%
\begin{equation*}
\mathcal{T}_{n-r}\mathcal{T}_{n+r}=\mathcal{T}_{2n}=\underset{n\text{ times}}%
{\underbrace{\mathcal{T}_{2}\mathcal{T}_{2}...\mathcal{T}_{2}}}=\mathcal{T}%
_{2}^{n}.
\end{equation*}
\end{description}

We have analogues results for the matrix sequence $\mathcal{K}_{n}.$

\begin{theorem}
For non-negatif integers $m,n$ and $r$ with $n\geq r,$ the following
identities hold:

\begin{description}
\item[(a)] $\mathcal{K}_{n-r}\mathcal{K}_{n+r}=\mathcal{K}_{n}^{2},$

\item[(b)] $\mathcal{K}_{n}^{m}=\mathcal{K}_{0}^{m}\mathcal{T}_{mn}.$
\end{description}
\end{theorem}

\textit{Proof.} \ \ 

\begin{description}
\item[(a)] We use Binet's formula of Tribonacci-Lucas matrix sequence which
is given in Theorem \ref{theorem:hbnvxgfdtreopas}. So%
\begin{eqnarray*}
&&\mathcal{K}_{n-r}\mathcal{K}_{n+r}-\mathcal{K}_{n}^{2} \\
&=&\mathcal{(}A_{2}\alpha ^{n-r}+B_{2}\beta ^{n-r}+C_{2}\gamma
^{n-r})(A_{2}\alpha ^{n+r}+B_{2}\beta ^{n+r}+C_{2}\gamma ^{n+r}) \\
&&-(A_{2}\alpha ^{n}+B_{2}\beta ^{n}+C_{2}\gamma ^{n})^{2} \\
&=&A_{2}B_{2}\alpha ^{n-r}\beta ^{n-r}(\alpha ^{r}-\beta
^{r})^{2}+A_{2}C_{2}\alpha ^{n-r}\gamma ^{n-r}(\alpha ^{r}-\gamma ^{r})^{2}
\\
&&+B_{2}C_{2}\beta ^{n-r}\gamma ^{n-r}(\beta ^{r}-\gamma ^{r})^{2} \\
&=&0
\end{eqnarray*}%
since $A_{2}B_{2}=A_{2}C_{2}=C_{2}B_{2}$ (see Lemma \ref%
{lemma:gbnhytrdvopxdser}). Now we get the result as required.

\item[(b)] By Theorem \ref{theorem:cvxhyugdswazxsa}, we have%
\begin{equation*}
\mathcal{K}_{0}^{m}\mathcal{T}_{mn}=\underset{m\text{ times}}{\underbrace{%
\mathcal{K}_{0}\mathcal{K}_{0}...\mathcal{K}_{0}}}\underset{m\text{ times}}{%
\underbrace{\mathcal{T}_{n}\mathcal{T}_{n}...\mathcal{T}_{n}}}.
\end{equation*}%
When we apply Lemma \ref{lemma:uygbcxfdzreawqops} (a) iteratively, it
follows that%
\begin{eqnarray*}
\mathcal{K}_{0}^{m}\mathcal{T}_{mn} &=&(\mathcal{K}_{0}\mathcal{T}_{n})(%
\mathcal{K}_{0}\mathcal{T}_{n})...(\mathcal{K}_{0}\mathcal{T}_{n}) \\
&=&\mathcal{K}_{n}\mathcal{K}_{n}...\mathcal{K}_{n}=\mathcal{K}_{n}^{m}.
\end{eqnarray*}%
This completes the proof.
\end{description}


\begin{thebibliography}{99}
\bibitem{bruce1984} \label{bib:bruce1984}Bruce, I., A modified Tribonacci
sequence, The Fibonacci Quarterly, 22 : 3, pp. 244--246, 1984.

\bibitem{cerdamoralez2018aa} \label{cerdamoralez2018aa}Cerda-Morales, G., On
the Third-Order Jabosthal and Third-Order Jabosthal-Lucas Sequences and
Their Matrix Representations, arxiv:1806.03709v1 [math.CO], 2018.

\bibitem{choi2013} \label{bib:choi2013}Choi, E., Modular tribonacci Numbers
by Matrix Method, J. Korean Soc. Math. Educ. Ser. B: Pure Appl. Math. Volume
20, Number 3 (August 2013), pages 207--221, 2013.

\bibitem{civciv2008} \label{civciv2008}Civciv, H., and Turkmen, R., On the
(s; t)-Fibonacci and Fibonacci matrix sequences, Ars Combin. 87, 161-173,
2008.

\bibitem{civciv2008b} \label{civciv2008b}Civciv, H., and Turkmen, R., Notes
on the (s; t)-Lucas and Lucas matrix sequences, Ars Combin. 89, 271-285,
2008.

\bibitem{feinberg1963} \label{bib:feinberg1963}Feinberg, M.,
Fibonacci--Tribonacci, The Fibonacci Quarterly, 1 : 3 (1963) pp. 71--74,
1963.

\bibitem{gulec2012} \label{gulec2012}Gulec, H.H., and Taskara, N., On the
(s; t)-Pell and (s; t)-Pell-Lucas sequences and their matrix
representations, Appl. Math. Lett. 25, 1554-1559, 2012.

\bibitem{howard2010} \label{bib:howard2010}Howard, F.T., Saidak, F.,
Congress Numer. 200 (2010), 225-237, 2010.

\bibitem{marohnic2012} \label{bib:marohnic2012}Marohni\'{c}, L., Strme\v{c}%
ki, T., Plastic Number: Construction and Applications, Advanced Research in
Scientific Areas 2012, 1523-1528, 2012.

\bibitem{piezas} \label{piezas}Piezas, T., A Tale of Four Constants,
https://sites.google.com/site/tpiezas/0012.

\bibitem{scott1977} \label{bib:scott1977}Scott, A., Delaney, T., Hoggatt
Jr.,\ V., The Tribonacci sequence, The Fibonacci Quarterly, 15:3, pp.
193--200, 1977.

\bibitem{shannon1977} \label{bib:shannon1977}Shannon, A., Tribonacci numbers
and Pascal's pyramid, The Fibonacci Quarterly, 15:3, pp. 268-275, 1977.

\bibitem{sloane} \label{bib:sloane}N.J.A. Sloane, The on-line encyclopedia
of integer sequences, http://oeis.org/

\bibitem{spickerman1981} \label{bib:spickerman1981}Spickerman, W., Binet's
formula for the Tribonacci sequence, The Fibonacci Quarterly, 20,
pp.118--120, 1981.

\bibitem{uslu2013} \label{uslu2013}Uslu, K., and Uygun, S., On the (s,t)
Jacobsthal and (s,t) Jacobsthal-Lucas Matrix Sequences, Ars Combin. 108,
13-22, 2013.

\bibitem{uygun2015} \label{uygun2015}Uygun, \c{S}., and Uslu, K.,
(s,t)-Generalized Jacobsthal Matrix Sequences, Springer Proceedings in
Mathematics\&Statistics, Computational Analysis, Amat, Ankara, May 2015,
325-336.

\bibitem{uygun2016} \label{uygun2016}Uygun, \c{S}., Some Sum Formulas of
(s,t)-Jacobsthal and (s,t)-Jacobsthal Lucas Matrix Sequences, Applied
Mathematics, 7, 61-69, 2016.

\bibitem{yalavigi1972} \label{bib:yalavigi1972}Yalavigi, C. C., Properties
of Tribonacci numbers, The Fibonacci Quarterly, 10 : 3, pp. 231--246, 1972.

\bibitem{yazlik2013a} \label{yazlik2013a}Yazlik, Y., and Taskara, N., Uslu,
K. and Yilmaz, N. The generalized (s; t)-sequence and its matrix sequence,
Am. Inst. Phys. (AIP) Conf. Proc. 1389, 381-384, 2012.

\bibitem{yilmaz2013} \label{yilmaz2013}Yilmaz, N., and Taskara, N., Matrix
Sequences in Terms of Padovan and Perrin Numbers, Journal of Applied
Mathematics, Volume 2013, Article ID 941673, 7 pages, 2013.

\bibitem{yilmaz2014a} \label{yilmaz2014a}Yilmaz, N., Taskara, N., On the
Negatively Subscripted Padovan and Perrin Matrix Sequences, Communications
in Mathematics and Applications, Vol. 5, No. 2, 59-72, 2014.

\bibitem{wani2018} \label{wani2018}Wani, A.A., Badshah, V.H., and Rathore,
G.B.S., Generalized Fibonacci and k-Pell Matrix Sequences, Punjab University
Journal of Mathematics (ISSN 1016-2526) Vol. 50(1) (2018) pp. 68-79.
\end{thebibliography}
\end{document}